\documentclass[12pt]{article}

\usepackage{amscd,amsmath, amssymb, fancyhdr,
	epsfig,color,bbm,url}
\definecolor{dblue}{rgb}{0,0,.6}

\usepackage[small,nohug,heads=littlevee]{diagrams}
\diagramstyle[labelstyle=\scriptstyle]

\usepackage[backref=page]{hyperref}
\renewcommand*{\backref}[1]{}
\renewcommand*{\backrefalt}[4]{%
	\ifcase #1 (Not cited.)%
	\or        (Cited on page~#2.)%
	\else      (Cited on pages~#2.)%
	\fi}

\hypersetup{
	colorlinks   = true,
	citecolor    = magenta,
	linkcolor= blue
}

\numberwithin{equation}{section}

\newcommand{\version}{version 1.0\ \   Sep. 8, 2019}

\def\eqref#1{(\ref{#1})}

\newcommand{\arrow}{{\:\longrightarrow\:}}
\newcommand{\Z}{{\Bbb Z}}
\newcommand{\C}{{\Bbb C}}

\newcommand{\R}{{\Bbb R}}

\def\1{\sqrt{-1}\:}

\newcommand{\cntrct}                
{\hspace{2pt}\raisebox{1pt}{\text{$\lrcorner$}}\hspace{2pt}}

\makeatletter
\def\x@arrow{\DOTSB\Relbar}
\def\xlongequalsignfill@{\arrowfill@\x@arrow\Relbar\x@arrow}
\newcommand{\xlongequal}[2][]{%
	\ext@arrow 0099\xlongequalsignfill@{#1}{#2}}
\def\xlongrightarrowfill@{\arrowfill@\relbar\relbar\longrightarrow}
\newcommand{\xlongrightarrow}[2][]{%
	\ext@arrow 0099\xlongrightarrowfill@{#1}{#2}}
\makeatother

\renewcommand{\phi}{\varphi}
\renewcommand{\epsilon}{\varepsilon}

\renewcommand{\leq}{\leqslant}

\newcommand{\im}{\operatorname{im}}

\newcommand{\diag}{\operatorname{\sf diag}}

\newcounter{Mycounter}[section]
\newcounter{lemma}[section]
\newcounter{claim}[section]
\newcounter{sublemma}[section]
\newcounter{corollary}[section]
\newcounter{theorem}[section]
\newcounter{conjecture}[section]
\newcounter{proposition}[section]
\newcounter{definition}[section]
\newcounter{example}[section]
\newcounter{remark}[section]
\newcounter{problem}[section]
\newcounter{question}[section]
\makeatletter

\renewcommand{\leftmark}{{\scriptsize  
		Dolbeault cohomology of compact complex manifolds}}

\@addtoreset{equation}{section} \@addtoreset{footnote}{section} \makeatother

\begin{document}
	\begin{center}
		{\LARGE\bf
			Dolbeault cohomology \\[3mm] of compact complex manifolds\\[4mm] with an action of a complex Lie \\[4mm] group.}
		\medskip
		\medskip

		Nikita Klemyatin
	\end{center}
	
	{\small \hspace{0.02\linewidth}
		\begin{minipage}[t]{0.85\linewidth} \small
			{\bf Abstract:}
			Let $G$ be a complex Lie group acting on a compact complex Hermitian manifold $M$ by holomorphic isometries. We prove that the induced action on the Dolbeault cohomology and on the Bott-Chern cohomology is trivial. We also apply this result to compute the Dolbeault cohomology of Vaisman manifolds.
		\end{minipage}
	}
	
	{\small 
		\tableofcontents
	}
	
	\section{Introduction.}
	
	One of the main invariants of a compact complex manifold is its Dolbeault cohomology groups $H^{*,*}_{\overline{\partial}}(M)$. They are defined as cohomology groups of the complex $(\Omega^{*,*}(M), \overline{\partial})$. Whenever $M$ is compact K\"ahler, Dolbeault cohomology are isomorphic to de Rham cohomology groups $H^*_{dR}(M; \C)$. More precisely, there is a Hodge decomposition for de Rham cohomology groups:
	
	$$H^k_{dR}(M; \C) = \bigoplus_{p+q=k} H^{p,q}_{\overline{\partial}}(M). $$
	
	However, when $M$ is non-K\"ahler, this decomposition generally fails.
	
	Now let $M$ be a manifold which is equipped with a smooth action of a connected Lie group $G$. A simple computation with Cartan formula $L_X = d \iota_X + \iota_X d$ shows that $G$ acts trivially on $H^*_{dR}(M)$. Hence, in the case of a compact K\"ahler $M$ with holomorphic action of $G$ the induced action of on $H^{*,*}_{\overline{\partial}}(M)$ is again trivial. In contrary, the non-K\"ahler case provides examples of non-trivial action (see \cite{_Akhiezer:Frst_}, \cite{_Akhiezer:Scnd_}) on the Dolbeaut cohomology groups. So one could ask when the holomorphic action of a Lie group induces the trivial action on $H^{*,*}_{\overline{\partial}}(M)$.
	
	This question could be extended on the {\it Bott-Chern} $H^{*,*}_{BC}(M)$ and {\it Aeppli} $H^{*,*}_{BC}(M)$ cohomology groups of a compact complex manifold $M$ with holomorphic action of a group $G$. These groups are defined as follows: $$H^{*,*}_{BC}(M) = \frac{\ker(d) \cap \ker(d^c)}{\im(dd^c)}$$ and $$H^{*,*}_A(M) = \frac{\ker(dd^c)}{\im(d) + \im(d^c)} .$$ 
	
	The induced action of $G$ on these cohomology groups could be non-trivial.
	However, there are conditions on $M$ and $G$, which implies the triviality of the induced action on $H^{*,*}_{\overline{\partial}}(M), H^{*,*}_{BC}(M)$ and   $H^{*,*}_A(M)$.

	\hfill
	
	\begin{theorem}\label{Main_theorem}
		Let $G$ be a complex Lie group, which acts on a compact Hermitian manifold $M = (M,h)$ by holomorphic isometries. Then $G$ acts trivially on Dolbeault, Bott-Chern and Aeppli cohomologies of $M$.
	\end{theorem}{}
	
	\hfill
	
	As a corollary, we obtain that the Dolbeault cohomology may be represented by invariant forms and any $\Delta_{\overline{\partial}}$-harmonic form is invariant.
	
	As an application of \ref{Main_theorem} we compute the Dolbeault cohomology of compact Vaisman manifolds. Indeed, Dolbeault cohomology groups of $M$ are isomorphic to the Dolbeault cohomology groups of invariant forms and they can be easily computed by this observation.
	
	\hfill
	
		\begin{theorem}\label{CohomologyOfVaisman}
		The Dolbeault cohomology groups of a Vaisman manifold $M$ are organized as follows:
		
		$$H_{\overline{\partial}}^{p,q}(M) = \begin{cases} \frac{H^{p,q}_{B}(M) \oplus \theta^{0,1}\wedge H^{p,q-1}_{B}(M)}{\text{Im}(L_{\omega_0})} ,~ p+q\leq \dim_\C(M) \\  \text{Ker}(L_{\omega_0})\vline_{H^{p,q}_{B}(M) \oplus \theta^{0,1}\wedge H^{p,q-1}_{B}(M)}, ~p+q>\dim_\C(M)\end{cases}$$
	\end{theorem}
	
	\section{Preliminaries.}
	\subsection{Basic properties of complex Lie groups.}
	We briefly introduce some notions about holomorphic group actions on complex manifolds.
	
	\begin{definition}
		A group $G$ is called a {\it complex Lie group} if $G$ is a complex manifold and the maps $G \times G \rightarrow G, (g_1; g_2) \mapsto g_1g_2$ and $G \rightarrow G, g \mapsto g^{-1}$ are holomorphic. 
	\end{definition}
	
	A complex Lie group $G$ is {\it compact} if the underlying complex manifold is compact. 
	
	The following claim is well-known.
	\hfill
	
	\begin{claim}\label{Torus-theorem}
		Let $G$ be a compact complex Lie group of dimension $n$. Then $G$ is a torus $\C^n/\Z^{2n}$.
	\end{claim}{}
	
	\hfill
	
	\begin{proof}
		Consider the adjoint representation $\text{Ad}: G \rightarrow \text{End}(\text{Lie}(G)) \cong \C^{n^2}$. This is a holomorphic map, hence it is constant by the maximum principle. Since $\text{Ad}(e) = \mathbbm{1}_n$, we have $\text{Ad}(G) = \mathbbm{1}_n$ (here $\mathbbm{1}_n$ is a $n \times n$ identity matrix). Hence $G$ is a commutative group and the exponential map $\exp: \text{Lie}(G) \rightarrow G$ is a homomorphism. Since $G$ is compact, $\exp$ has a nontrivial kernel, which is a finitely generated abelian group without torsion, i.e. $\ker(\exp) \cong \Z^p$. It is easy to see, that $p=2n$ and $G \cong \C^n/\Z^{2n}$.
	\end{proof}{}
	
	\hfill
	
	When $G$ is non-compact, it might be non-abelian.
	However, when $G$ acts on $M$ by holomorphic isometries, we can say something about the closure of $G$. 
	First of all we need the following result about the isometry group of a compact Riemannian manifold.
	
	\hfill
	
	\begin{theorem}\label{IsometryGroupOfManifold}(see \cite{Helgason}, Chapter II, Thm.1.2)
		Let $(M,g)$ be a compact Riemannian manifold. Then the group  $\text{Isom}(M)$ isometries of $(M,g)$ is a compact Lie group.
	\end{theorem}
	
	\hfill
	
	We also need a following statement.
	
	\hfill
	
	\begin{proposition}\label{ClosureOfGroup}
		Let $G$ be a complex group which acts by a holomorphic isometries on a compact Hermitian manifold $(M;h)$. Then the closure of $G$ in the group $\text{Isom}(M)$ still acts by holomorphic isometries on $M$.
	\end{proposition}{}
	
	\hfill
	
	\begin{proof}
		Let $\{g_n\}^{\infty}_{n=1}$ be a Cauchy sequence in $G$. Since each $g_n$ lies in $\text{Isom}$, the  family $\{g_n\}^{\infty}_{n=1}$ is bounded.
		By the Montel's theorem, the limit $g$ of $\{g_n\}^{\infty}_{n=1}$ is holomorphic. By the theorem of Cartan and von Neumann (see \cite{Helgason}), the closure is a Lie group. It will be compact by \ref{IsometryGroupOfManifold}.
	\end{proof}{}
	
	\subsection{Holomorphic vector fields and Lie derivative.}
	
	Let $G$ be a Lie group and $M = (M,J)$ be a compact complex manifold. We assume that both $G$ and $M$ are connected.
	
	\hfill
	
	\begin{definition}
		A holomorphic action of complex Lie group $G$ on a complex manifold $M$ is a holomorphic map $\sigma: G \times M \rightarrow M$.
	\end{definition}
	
	\hfill
	
	\begin{definition}(\cite{_Gauduchon:Book_})
		A vector field $X$ is called {\it holomorphic vector field} if $L_X J = 0$.
	\end{definition}{}
	
	\hfill
	
	Clearly holomorphic vector fields form a Lie subalgebra $\mathfrak{h}(M)$ in the algebra of all holomorphic vector fields on $M$. If $G$ acts holomorphically on $M$, then the action induces a homomorphism $d_e\sigma: \text{Lie}(G) \rightarrow \mathfrak{h}(M)$ of (real) Lie algebras.
	This homomorphism gives us the action of $\text{Lie}(G)$ on tensors on $M$ via the Lie derivative. Recall that Lie derivative acts on differential forms via {\it Cartan formula}:
	
	$$L_X \alpha= d\iota_X \alpha + \iota_X d \alpha.$$
	
	It also commutes with de Rham differential: $$ L_Xd \alpha = dL_X \alpha.$$
	
	Also the Lie derivative satisfies the Leibniz rule for the tensor product:
	
	$$L_X(P \otimes Q) = (L_X P) \otimes O + P \otimes (L_X Q).$$
	
	In particulary, the following identity holds: $$ L_X Y = [X;Y].$$
	
	As mentioned above, all holomorphic vector fields form a Lie algebra. Moreover, on a complex manifold, the algebra of holomorphic vector fields is a {\it complex} Lie algebra and $JX$ is holomorphic whenever $X$ is holomorphic. Indeed, for a holomorphic vector field $X$ and for any vector field $Y$ one has 
	\begin{align*}
	L_{JX}(JY) = [JX;JY] = J[JX;Y] + J[X;JY] + [X;Y] = \\
	= J[JX;Y] - [X;Y] + [X;Y] = J[JX;Y] = JL_XY.
	\end{align*}{}
	Here in the second equality we use that Nijenhuis tensor is zero and in the third equality we use the definition of holomorphic vector fields.
	
	Now, in the case of $G$-action on $M$ that the homomorphism $d_e\sigma: \text{Lie}(T) \rightarrow \mathfrak{h}(M)$ is actually a homomorphism of {\it complex} Lie algebras. 
	
	\subsection{$d^c$-operator and Lie derivative.}
	
	We start from the following definition:
	
	\hfill
	
	\begin{definition}(see \cite{_Gauduchon:Book_}, sect. 1.11)
		Consider an operator $d^c$, defined as following:
		
		$$d^c \alpha = JdJ^{-1}\alpha = (-1)^{|\alpha |} JdJ\alpha. $$
		Here $\alpha$ is an arbitrary form of degree $| \alpha |$.
	\end{definition}
	
	\hfill
	
	Note that $d^c$ is a real operator. 
	
	One can write down $d^c$ in terms of operators $\partial$ and $\overline{\partial}$ (see \cite{_Gauduchon:Book_}): $$d^c = i(\overline{\partial} - \partial)$$. Here $d=\partial + \overline{\partial}$ is the standard decomposition of the de Rham differential on complex manifolds.
	
	We also have 
	
	\begin{align*}
	\overline{\partial} = \frac{1}{2}(d + id^c)
	\hspace{10 mm}
	\partial = \frac{1}{2}(d - id^c).
	\end{align*}{}
	
	There is a simple Cartan-type formula, which connects the operator $d^c$ and the 
	Lie derivative along $L_{JX}$.
	
	\hfill
	
	\begin{proposition}\label{aLaCartan}
		Let $X$ be a holomorphic vector field on compact complex manifold $M$. Then $L_{JX}\alpha = -\{d^c;\iota_X\} = -(d^c \iota_X + \iota_X d^c)\alpha$
		for any form $\alpha$. 
	\end{proposition}{}
	
	\hfill
	
	\begin{proof}
		$$L_{JX} \alpha = (-1)^{|\alpha|}L_{JX}J^2 \alpha = (-1)^{|\alpha|}J\{d; \iota_{JX}\}J \alpha = (-1)^{|\alpha| + 1}\{JdJ; \iota_{X}\}\alpha = $$

		$$ = -\{d^c; \iota_{JX}\}\alpha. $$
		
	\end{proof}
	
	\subsection{Elliptic operators on compact manifolds.}
	
	In this section we recall some facts about elliptic operators on compact manifolds. 
	
	Let $M$ be a compact manifold and $E$ a vector bundle on $M$. Throughout this section the symbol $\Gamma(E)$ denotes the space of smooth sections on $M$.
	
	\hfill
	
	\begin{definition}
		Let $E,F$ are complex vector bundles over compact manifold $M$. A {\it linear differential operator} $P$ of degree $k$ from $E$ to $F$ is a $\C$-linear operator $P: \Gamma(E) \rightarrow \Gamma(F), ~s \mapsto Ps$ of the form
	
	\begin{equation}
		Ps(x) = \sum_{|j|=0}^k a_j(x) \frac{\partial^j}{\partial x_1^{j_1} \dots \partial x_m^{j_m}}s(x).
	\end{equation}
	
	\end{definition}
	
	The {\it principal symbol} of the operator $P$ is the morphism of vector bundles $\sigma_P(x, \xi) = \sum_{|j|=k} a_j(x) \xi^j, ~\sigma_P(x, \xi): E \rightarrow F$. Here $\xi \in T_\xi M$.  
	
	The operator $P$ is called {\it elliptic}, if the principal symbol $\sigma_P(x, \xi)$ of $P$ is injective for any nonzero $\xi$.
	
	Fix a volume form $d\mu$ on $M$. Let $E$ be a Hermitian vector bundle,  the operator $P: \Gamma(E) \rightarrow \Gamma(E)$ is called self-adjoint, if it is self-adjoint with respect to the $L^2$ scalar product 
	
	$$(s;t) := \int_M h(s,t) d\mu$$
	
	for any two $s,t \in \Gamma(E)$.

	Elliptic self-adjoint operators have very nice spectral properties.
	
	\hfill
	
	\begin{proposition}\label{SpectralDecomposition}(see \cite{_Gilkey:Book_})
		Let $P: \Gamma(E) \rightarrow \Gamma(E)$ be an elliptic self-adjoint operator. Then we can find a complete basis $\{s_j\}_{j=1}^{\infty}$ of $L_2(E)$ of eigensections of $P$.  Each eigensection of $P$ is smooth and each eigenspace has finite dimension. Moreover the set of eigenvalues is a discrete subset in $\R$.
	\end{proposition}

	\hfill

	This proposition has a very nice application for the case of $E= \Omega^*(M)$ on a compact Hermitian manifold $M$ with the holomorphic and isometric action of compact Lie group $G$.
	
	\hfill	
	
	\begin{proposition}\label{RepsOnEigenspaces}
		Let $G$ and $M$ are as above and $P: \Omega^*(M) \rightarrow \Omega^*(M)$ be a self-adjoint elliptic operator. Suppose that for any $g \in G$ and for any smooth form $\alpha$ we have $g^* P\alpha = Pg^* \alpha$. Then each eigenspace of $P$ is a direct sum of irreducible representations of $G$. Moreover, for any $g \in G$ there exist a complete orthogonal basis $\{\alpha_j\}_{j=0}^{\infty}$ such that $P\alpha_j = \lambda_j \alpha_j$ and $g^* \alpha_j = e^{ia_j}\alpha_j, ~a_j \in \R$.
	\end{proposition}

	\hfill
	
	\begin{proof}
		Since $g^* P\alpha = Pg^* \alpha$, $g^*$ and $P$ commute on each eigenspace of $P$. Hence they have common eigenvectors on each eigenspace. For an arbitrary $g \in G$ the map $g^*$ is $L^2$-isometry and hence its restriction on each eigenspace is an unitary operator. Hence all eigenvalues of $g^*$ are equal to $e^{ia_j}$ for some real $a_j$.
	\end{proof}

	\subsection{Various cohomology groups on compact complex manifolds.}
	
	Recall that Dolbeault cohomology groups are defined as follows: 
	
	$$H^{*,*}_{\overline{\partial}}(M) = \frac{\ker(\overline{\partial})}{\im(\overline{\partial})}.$$
	
	In the case of K\"ahler manifolds there is a decomposition of de Rham cohomology groups on a direct sum of Dolbeault cohomology groups:
	
	$$H^k_{dR}(M; \C) = \bigoplus_{p+q=k} H^{p,q}_{\overline{\partial}}(M). $$
	
	In general, there is no such decomposition for de Rham cohomology. However, there is the Frolicher inequality which is obtained from the Fr\"olicher spectral sequence (see \cite{_Moroianu:Book_}): 
	
	$$\dim_{\C}H^k_{dR}(M; \C) \leq \sum_{p+q=k}\dim_{\C}H^{p,q}_{\overline{\partial}}(M). $$
	
	One can define another cohomology groups, namely {\it Bott-Chern} 
	cohomology $$H^{p,q}_{BC}(M) = \frac{\ker(d) \cap \ker(d^c)}{\im(dd^c)} $$ and 
	{\it Aeppli} cohomology $$H^{p,q}_A(M) = \frac{\ker(dd^c)}{\im(d) + \im(d^c)} .$$
	
	There is an analogue of harmonic decomposition for these cohomologies. Indeed, for any Hermitian metric $h$ on $M$ one can construct the following Laplacians for Bott-Chern and Aeppli cohomologies (see \cite{_Schweitzer_}):
	
	$$\Delta_{BC} := (\partial \overline{\partial})(\partial \overline{\partial})^* + (\partial \overline{\partial})^*(\partial \overline{\partial}) + (\overline{\partial}^* \partial)(\overline{\partial}^* \partial)^* + (\overline{\partial}^* \partial)^*(\overline{\partial}^* \partial) + \overline{\partial}^* \overline{\partial} + \partial^* \partial $$
	
	and
	
	$$\Delta_A = \overline{\partial} \overline{\partial}^* + \partial \partial^* +  (\partial \overline{\partial})(\partial \overline{\partial})^* + (\partial \overline{\partial})^*(\partial \overline{\partial}) + (\overline{\partial} \partial^*)(\overline{\partial} \partial^*)^* + (\overline{\partial} \partial^*)^*(\overline{\partial} \partial^*).$$
	
	These operators are self-adjoint and elliptic and their kernels are isomorphic to $H^{*,*}_{BC}(M)$ and $H^{*,*}_{A}(M)$ respectively.
	
	Both Bott-Chern and Aeppli cohomologies are ``dual'' to each other in the following sense.
	
	\hfill
	
	\begin{theorem}(see \cite{_Angela_}, Theorem 2.5, and \cite{_Schweitzer_})\label{_Nondegenerate_}
		Let $M$ be a compact complex manifold of complex dimension $m$. Then there is an isomorphism between $H^{p,q}_{BC}(M)$ and $(H^{m-p,m-q}_{A}(M))^*$. The isomorphism is given by the nondegenerate pairing $$\int_M : H^{p,q}_{BC}(M) \times H^{m-p,m-q}_{A}(M) \rightarrow \C. $$ 
	\end{theorem}
	
	\hfill
	
	All these cohomology groups are related in the following way (see \cite{_AngTom_}):
	
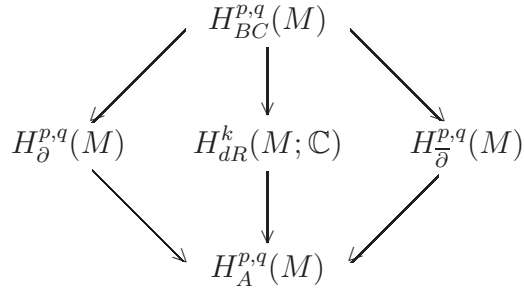
	\begin{diagram}
		& & & & H^{p,q}_{BC}(M) & & & &\\
		& & \ldTo &     &    \dTo    &     &  \rdTo & &\\
		&H^{p,q}_{\partial}(M) & & & H^k_{dR}(M; \C)  & &  &  H^{p,q}_{\overline{\partial}}(M)  & \\
		& & \rdTo &     &    \dTo    &     &  \ldTo & &\\
		& & & & H^{p,q}_{A}(M) & & & &\\
	\end{diagram}
	
	Note that in the general case all arrows in this diagram are neither injective nor surjective. For instance, if $M$ is non-K\"ahler, then the map $H^{p,q}_{BC}(M) \rightarrow H^k_{dR}(M)$ may have a nontrivial kernel (see \cite{_OVV_}, Theorem 2.3).

	\section{Proof of the main theorem}
	
	\begin{proposition} \label{_Bott-Chern_}
		Suppose that complex connected Lie group $G$ acts by holomorphic isometries on compact complex non-K\"ahler Hermitian manifold $(M,h)$. Then $G$ acts trivially on $H^{*,*}_{BC}(M)$.
	\end{proposition}
	
	\hfill
	
	\begin{proof}
		Consider the natural map $H^{*,*}_{BC}(M) \rightarrow H^{*}_{dR}(M; \C)$. Since $G$ acts trivially on $H^{*}_{dR}(M;\C)$, it is sufficient to prove that the action on $\text{Ker}(H^{*,*}_{BC}(M) \rightarrow H^{*}_{dR}(M; \C))$ is trivial.
		
		Suppose now that $\eta = d\alpha$ for some $\alpha$. Let $X \in \text{Lie}(G)$. Each $\alpha$ can be written as $\alpha = \sum_{j=0}^{\infty}c_j \alpha_j$, there $\Delta_{BC}\alpha_j = \lambda_j\alpha_j$. We also have $L_X\alpha = \sum_{j=0}^{\infty}c_j L_X\alpha_j$. Hence, by \ref{RepsOnEigenspaces} we can assume $L_X \alpha = a\alpha$ and $L_{JX}\alpha = b\alpha$ for some $a,b \in i\R$.
		
		Since $d^c \eta = 0$, we have
		\begin{align*}
		0 = \iota_X d^c d \alpha = -L_{JX}d\alpha - d^c \iota_X d \alpha = -L_{JX}d\alpha - \\-d^c L_X \alpha + d^cd\iota_X \alpha  =  -(bd + ad^c) \alpha + d^cd\iota_X \alpha.
		\end{align*}
		
		Let $\delta := bd+ad^c$, where $a,b \in i\R$ are as above. From the equation above we can see that $0 = -\delta \alpha + d^cd \iota_X \alpha$. 
		
		Now we want to write down $L_X$ in the form $\{\delta, \iota_Y\}$ for some $Y= y_1 X + y_2 JX$:
		\begin{align*}
		\{\delta, \iota_Y\} = by_1\{d^c, \iota_X\} + ay_1\{d;\iota_X\} + ay_2\{\iota_{JX};d\} + by_2\{\iota_{JX};d^c\} = \\ = (ay_2 - by_1)L_{JX} + (ay_1 + by_2)L_X = L_X.
		\end{align*}
		
		Therefore we have the system of linear equations:
		\[
		\begin{cases} 
		ay_2 - by_1 = 0 \\
		ay_1 + by_2 = 1. \end{cases}
		\]
		The determinant of this system is equal to $-b^2 - a^2$ and it is nonzero whenever either $a$ or $b$ is nonzero. Hence this system has a solution, i.e. $L_X \alpha = \{\delta, \iota_Y\} = \delta \iota_Y \alpha$.
		
		If $\eta = d\alpha$ then $L_X \eta = L_X d \alpha = d L_X \alpha = d \delta \iota_Y \alpha = add^c\iota_Y \alpha = dd^c \beta$ for some $\beta$. Thus $L_X$ does not change the class of $\eta$ in $H^{*,*}_{BC}(M)$.
		
	\end{proof}
	
	\hfill
	
	\begin{corollary}\label{_Aeppli_}
		In the assumptions of \ref{_Bott-Chern_}, $G$ acts trivially on $H^{*,*}_{A}(M)$.
	\end{corollary}
	
	\hfill
	
	\begin{proof}
		This is a direct corollary of \ref{_Bott-Chern_} and \ref{_Nondegenerate_}.
	\end{proof}
	
	\hfill
	
	Now we can prove \ref{Main_theorem}.
	
	\hfill
	
	\begin{theorem}(also \ref{Main_theorem})
		Let $G$ be a compact Lie group which acts on a compact Hermitian manifold $M = (M,h)$ by holomorphic isometries. Then $G$ acts trivially on Dolbeault, Bott-Chern and Aeppli cohomologies of $M$.
	\end{theorem}{}
	
	\hfill
	
	\begin{proof}
		
		Consider the exact sequence:
		\begin{diagram}
			& 0 & \arrow & A^{p,q} & \arrow & B^{p,q} & \arrow & H^{p,q}_{\overline{\partial}}(M) & \arrow & H^{p,q}_{A}(M) & \arrow & C^{p,q} & \arrow & 0 &
		\end{diagram}
		(see \cite{_AngTom_}).
		
		Here we define $A^{p,q}, B^{p,q}$ and $C^{p,q}$ as follows:
		
		\begin{align*}
		A^{p,q} = \frac{\im(\partial) \cap \im(\overline{\partial})}{\im(\partial\overline{\partial})}\hspace{10 mm}
		B^{p,q} = \frac{\im(\partial) \cap \ker(\overline{\partial})}{\im(\partial\overline{\partial})}
		\hspace{10 mm}
		C^{p,q} = \frac{\ker(\partial\overline{\partial})}{\im(\partial)+\im(\overline{\partial})}.
		\end{align*}{}
		
		Since $A^{p,q}$ and $B^{p,q}$ are subspaces of $H^{p,q}_{BC}(M)$, the action of $G$ on $A^{p,q}$ and $B^{p,q}$ is trivial. The $G$-action on groups $H^{p,q}_{A}(M)$ are also trivial by \ref{_Aeppli_}. The groups $C^{p,q}$ have trivial $G$-action because the map $H^{p,q}_{A}(M) \rightarrow C^{p,q}$ is surjective and commutes with action of $G$. Hence the induced $G$-action on each group $H^{p,q}_{\overline{\partial}}(M)$ is trivial.

	\end{proof}
	
	\hfill
	
	\begin{corollary}\label{HarmonicInvariance} Let $h$ be a $G$-invariant Hermitian metric on $M$ and $\alpha \in H^{p,q}_{\overline{\partial}}(M)$ be a harmonic form (with respect to $h$) on $M$. Then $\alpha$ is $G$-invariant.
	\end{corollary}
	
	\hfill
	
	\begin{proof}
		Let $X$ be a holomorphic vector field from $\text{Lie}(G)$. Since
		$h$ is $G$-invariant, $L_X$ commutes with
		$\Delta_{\overline{\partial}}$ and hence $G$ acts on
		$\Delta_{\overline{\partial}}$-harmonic forms. However, the space
		of $\Delta_{\overline{\partial}}$-harmonic forms is
		isomorphic to the Dolbeault cohomology (see
		\cite{_Moroianu:Book_}). Since $g^*\alpha$ is again
		$\Delta_{\overline{\partial}}$-harmonic and has the same
		Dolbeault cohomology class as $\alpha$, we have $g^*
		\alpha = \alpha$.
	\end{proof}{}

	\hfill
	
	\begin{corollary}\label{InvariantForms}
		Let $M$ and $G$ be as above. Then the closure $\overline{G}$ of $G$ in $\text{Isom}(M)$ acts trivially on $H^{p,q}_{\overline{\partial}}(M)$.
	\end{corollary}

	\hfill
	
	\begin{proof}
		By \ref{ClosureOfGroup} the group $\overline{G}$ acts on the space of $\Delta_{\overline{\partial}}$-harmonic forms. The group $\overline{G}$ acts trivially on harmonic forms because $G$ acts trivially on it. 
	\end{proof}
	
	\section{The application to Vaisman manifolds.}
	
	\subsection{Sasakian and Vaisman manifolds.}
	
	In this section we recall some facts about Vaisman and Sasakian manifolds.
	
	\hfill
	
	\begin{definition}\label{SasDef}
		An odd-dimensional Riemann manifold $(S,g)$ is called a {\it Sasakian} manifold if it cone $C(S) :=S \times \R_{>0}$ with the metric $\widetilde{g}:=t^{2}g + dt^2$ is K\"ahler and the natural action of $\R^{>0}$ on $C(S)$ is holomorphic.
	\end{definition}{}
	
	\hfill
	
	\begin{definition}\label{LCKDef}
		A compact complex Hermitian manifold $(M, J, \omega)$ of $\text{dim}_\C >1$  is called a {\it locally conformally K\"ahler} (LCK for short), if it admits a K\"ahler covering $(\widetilde{M}, \widetilde{J}, \widetilde{\omega})$, such that covering group acts by holomorphic homoteties on $\widetilde{M}$.
	\end{definition}{}
	
	\hfill
	
	The LCK property is equivalent to existence of a closed form $\theta$ such that $d \omega = \theta \wedge \omega$. The form $\theta$ is called the {\it Lee form}. It is obviosly closed. When $\theta$ is exact, an LCK manifold can be equipped a K\"ahler metric. Indeed, if $\theta = d \phi$, then $e^{-\phi}\omega$ is closed.
	
	A very important example of LCK manifolds are Vaisman manifolds.
	
	\hfill
	
	\begin{definition}\label{VaisDef}
		A {\it Vaisman manifold} $(M, J, \omega, \theta)$ is an LCK manifold such that the Lee form $\theta$ is parallel with respect to the Levi-Civita connection which is associated to the Hermitian metric.
	\end{definition}{}
	
	\hfill
	
	The typical examples of Vaisman manifolds are Hopf varieties $H_A := (\C^n \setminus {0}) / \langle A \rangle$, where $A=\diag(\lambda_i)$ with $|\lambda_i|<1$ (see \cite{_OrneaVerbit:Second_}).
	
	The Vaisman manifold has a foliation $\Sigma$ which is called the {\it canonical (or fundamental) foliation}. It is generated by $X = \theta^\sharp := g^{-1}(\theta, \cdot)$ and $JX = J\theta^\sharp$. It is well-known that $X$ and $JX$ acts holomorphically on $M$. Moreover, there is a {\it transversely K\"ahler} metric on $M$. It is given by the following formula:
	
	$$2\omega_0 = d\theta^c = d(J\theta) = \omega - \theta \wedge \theta^c . $$
	
	 (see \cite{_Verbit:First_}). The local structure of compact Vaisman manifolds is well-known and it is described by the following theorem.
	
	\hfill
	
	\begin{theorem}\label{VaismanStructure}(The local structure theorem for compact Vaisman manifolds, see \cite{_OrneaVerbit:First_}) Let $(M, J, \omega, \theta)$ be a Vaisman manifold. Denote by $X$ the vector field dual to $\theta$. Then $L_XJ=0$ and $M$ is locally isomorphic to the K\"ahler cone of a Sasakian manifold. Moreover, $X$ acts on K\"ahler covering $\widetilde{M}$ by holomorphic homotheties of K\"ahler metric.
	\end{theorem}{}
	
	\hfill

	The following proposition is well-known.
	
	\hfill
	
	\begin{proposition}\label{Sas}
		The vector field $JX$ is a Killing.
	\end{proposition}

	\hfill
	
	\begin{proof}
		This is a local statement.
		
		By \ref{VaismanStructure} we can assume that locally $M=S \times \R$ with product metric $g_S + dt^2$. Moreover, we can assume that $\theta = dt$ and $X= \frac{d}{dt}$. Hence the metric $\widetilde{g} = e^{-t}(g_S + dt^2)$ is K\"ahler. Denote by $\widetilde{\omega}$ a K\"ahler form of $\widetilde{g}$.
		
		Since $JX$ is ortogonal to $X$, it is tangent to $S$. Since $JX$ is holomorphic, we have $L_{JX}\widetilde{g} = L_{JX} \widetilde{\omega} = d \iota_{JX} \widetilde{\omega} = d \theta = 0$.
	\end{proof}

	\hfill
	
	Since $\nabla X = \nabla \theta = 0$, $X$ is Killing as well and we have the following statement
	
	\hfill
	
	\begin{claim}\label{VaismanGroup}
		The group generated by $e^{tX}$ and $e^{tJX}$  acts by holomorphic isometries on $M$.
	\end{claim}
	\subsection{Dolbeault cohomology of Vaisman manifolds.}
	
	Let $(M, J, \omega, \theta)$ be a Vaisman manifold. We start from the corollary of \ref{Main_theorem}.
	
	\hfill
	
	\begin{corollary}
		The Dolbeault cohomology groups of $M$ are the cohomology groups of complex $(\Lambda^*(M)^{\text{inv}}, \overline{\partial})$ of invariant forms on $M$.
	\end{corollary}{}
	
	\hfill
	
	\begin{proof}
		The group, generated by $e^{tX}$ and $e^{tJX}$ acts by holomorphic isometries on $M$ (\ref{VaismanGroup}). Hence we can apply \ref{Main_theorem} and \ref{InvariantForms}. Each $\overline{\partial}$-closed invariant form lies in $\Lambda^*(M)^\Sigma$ and it is an element of cohomology group of the complex $(\Lambda^*(M)^{\text{inv}}, \overline{\partial})$. 
	\end{proof}{}

	\hfill
	
	Recall the important definition.
	
	\hfill
	
	\begin{definition}
		Let $M$ be a manifold with foliation $\Sigma$. A form $\alpha$ is {\it basic}, if $\iota_X \alpha = \iota_X d \alpha = 0$ for any vector field $X$ tangent to $\Sigma$.
	\end{definition}{}
	
	\hfill
	
	\begin{proposition}\label{InvForms}
		Let $\eta$ be an invariant form on $M$. Then $\Lambda^*(M)^{\text{inv}} = (\pi^*\Lambda^*_B)^{\text{inv}} \otimes \Lambda[\theta^{1,0}, \theta^{0,1}]$.
		
	\end{proposition}{}
	
	\hfill
	
	\begin{proof} A form $\alpha$ on $M$ is $G$-invariant iff it is invariant under the induced action of $\text{Lie}(G)$. Hence this is a purely local statement.
		Denote by $F$ the fiber of foliation $\Sigma$
		on $M$ and by $B$ the leaf space.
		
		We know that locally $$\Lambda^*(M) = \pi^*\Lambda^*_B \otimes \Lambda^*(F). $$
		
		Hence, the following equality holds for invariant forms
		
		$$\Lambda^*(M)^{\text{inv}} = (\pi^*\Lambda^*_B)^{\text{inv}} \otimes \Lambda^*(F)^{\text{inv}}. $$
		
		But the $(\pi^*\Lambda^*_B)^{\text{inv}}$ are just basic forms and $\Lambda^*(F)^{\text{inv}}$ is the exterior algebra generated by $\theta$ and $\theta^c$.		
	\end{proof}{}

	\hfill

	Recall some important definition from homological algebra.
	
	\hfill
	
	\begin{definition}\label{Cone}
		Suppose $(K^*,d_K)$ and $(L^*,d_L)$ are complexes and $f: K^* \rightarrow L^*$ be a morphism of these complexes. Define a complex $(C(f), d_f)$ as follows: $C(f)_i = K_{i+1} \oplus L_i$ and $ d_f = (d_K, f - d_L)$. This complex is called {\it the cone of} $f$.
	\end{definition}
	
	\hfill
	
	For each cone of a morphism we can construct the long exact sequence of cohomology. Indeed, we have the short exact sequence of complexes:
	
	$$0 \rightarrow L^* \rightarrow C(f) \rightarrow K^*[1] \rightarrow 0. $$
	
	There is a well-known way to construct a long sequence in cohomology from a short sequence of complexes:
	
	\begin{diagram}
		& \dots & \arrow & H^i(L^*) & \arrow & H^i(C(f)) & \arrow H^{i+1}(K^*) & \arrow H^{i+1}(K^*) & \arrow & \dots & 
	\end{diagram}

	(see \cite{_Gelfand:Manin_} for the details).
	
	\hfill
	
	Now we can compute the Dolbeault cohomology groups for a Vaisman manifold $M$.
	
	Consider the subcomplex $\Lambda^{*,*}_{B, \theta^{0,1}} := (\pi^*\Lambda^*_B)^{\text{inv}} \oplus \theta^{0,1} \wedge (\pi^*\Lambda^*_B)^{\text{inv}}$ of the complex $\Lambda^*(M)^{\text{inv}}$. Denote by $L_{\omega_0}$ the operator of multiplication by $\omega_0$. Clearly, $\Lambda^*(M)^{\text{inv}} = \Lambda^{*,*}_{B, \theta^{0,1}} \oplus \theta^{1,0} \wedge \Lambda^{*,*}_{B, \theta^{0,1}}$ and $L_{\omega_0}$ is a morphism $\Lambda^{*,*}_{B, \theta^{0,1}} \longrightarrow \Lambda^{*+1,*+1}_{B, \theta^{0,1}}$.
	
	\begin{proposition}\label{ConeInv}(\cite{_OrneaVerbit:Third_} in the case of de Rham cohomology)
		The complex $\Lambda^*(M)^{\text{inv}}$ is isomorphic to the cone $C(L_{\omega_0})$ of the morphism $\Lambda^{*,*}_{B, \theta^{0,1}} \xlongrightarrow{L_{\omega_0}} \Lambda^{*,*}_{B, \theta^{0,1}}$.
	\end{proposition}

	\hfill
	
	\begin{proof}
		We have $$\Lambda^*(M)^{\text{inv}} = \Lambda^{*,*}_{B, \theta^{0,1}} \oplus \theta^{1,0} \wedge \Lambda^{*,*}_{B, \theta^{0,1}} = \Lambda^{*,*}_{B, \theta^{0,1}} \oplus \Lambda^{*,*}_{B, \theta^{0,1}}[-1].$$
		
		The Dolbeault differential $\overline{\partial}$ on $\Lambda^*(M)^{\text{inv}}$ acts in the following way: it is the ordinary $\overline{\partial}$ on $\Lambda^{*,*}_{B, \theta^{0,1}}$. On the other hand, we have $\overline{\partial} \theta^{1,0} = \omega_0$. Hence the Dolbeault differential acts on $\Lambda^{*,*}_{B, \theta^{0,1}}[-1]$ as $\omega_0 - \overline{\partial}$.
	\end{proof}

	\hfill
	
	\begin{theorem}(also \ref{CohomologyOfVaisman})
		The Dolbeault cohomology groups of a Vaisman manifold $M$ are organized as follows:
		
		$$H_{\overline{\partial}}^{p,q}(M) = \begin{cases} \frac{H^{p,q}_{B}(M) \oplus \theta^{0,1}\wedge H^{p,q-1}_{B}(M)}{\text{Im}(L_{\omega_0})} ,~ p+q\leq \dim_\C(M) \\  \text{Ker}(L_{\omega_0})\vline_{H^{p,q}_{B}(M) \oplus \theta^{0,1}\wedge H^{p,q-1}_{B}(M)}, ~p+q>\dim_\C(M)\end{cases}$$
	\end{theorem}

	\hfill
	
	This result is similar to Theorem 3.2 from \cite{_Tsukada_}.
	
	\hfill
	
	\begin{proof} 
		We have a long exact sequence:
			\begin{diagram}
			& \dots & \rightarrow & H^{p,q}_{B, \theta^{0,1}}(M) & \xrightarrow{L_{\omega_0}} & H^{p+1,q+1}_{B, \theta^{0,1}}(M) & \rightarrow & H^{p+1,q+1}_{\overline{\partial}}(M) & \rightarrow & H^{p+1,q+1}_{B, \theta^{0,1}}(M) & \rightarrow & \dots & 
		\end{diagram}
		
		The cohomology groups $H^{*,*}_{B, \theta^{0,1}}(M)$ of complex $(\Lambda^{*,*}_{B, \theta^{0,1}}; \overline{\partial})$ are equal to $H^{*,*}_{B}(M) \oplus \theta^{0,1} \wedge H^{*,*-1}_{B}(M)$. Since $H^{*,*}_{B}(M)$ admits a Lefshetz $\text{SL}(2)$-action (see \cite{_ElKacimi_} and \cite{_ElKacimi:Hector}), there is an analog of such action for $H^{*,*}_{B, \theta^{0,1}}(M)$. Since $L_{\omega_0}$ is injective on $H^{p,q}_{B}(M)$ for $p+q \leq \dim_{\C}(M)$, it is injective on $H^{p,q}_{B, \theta^{0,1}}(M)$ with the same $p,q$. Hence we obtain the short exact sequence:
		\begin{diagram}
			& 0 & \arrow & H^{p-1,q-1}_{B, \theta^{0,1}}(M) & \xlongrightarrow{L_{\omega_0}} & H^{p,q}_{B, \theta^{0,1}}(M) & \arrow & H^{p,q}_{\overline{\partial}}(M) & \arrow & 0 & 
		\end{diagram}
	
		For the case $p+q > \dim_{\C}(M)$ we have another short exact sequence:
		\begin{diagram}
			& 0 & \arrow & H^{p,q}_{\overline{\partial}}(M) & \arrow   H^{p,q}_{B, \theta^{0,1}}(M) & \xlongrightarrow{L_{\omega_0}} & H^{p+1,q+1}_{B, \theta^{0,1}}(M) & \arrow & 0 & 
		\end{diagram}
	
	The statement of the theorem directly follows from these two sequences.
		
	\end{proof}

	\hfill

	{\bf Acknowledgements:} Many thanks to Misha Verbitsky for fruitful discussions and help. I also want to thank Instituto Nacional de Matemática Pura e Aplicada for hospitality.

	\footnotesize
	{
		
		\noindent
		{\sc Nikita Klemyatin\\
			\sc National Research University HSE,\\
			Department of Mathematics, 6 Usacheva Str. Moscow, Russia\\
			also:\\
			\sc Skolkovo Institute of Science and Technology\\
			Bolshoy Boulevard 30, bld. 1. Moscow, Russia\\			
			\tt  nklemyatin@yandex.ru}

	}
	
\end{document}